\documentclass[12pt]{article}
\usepackage{amssymb}
\usepackage{amsmath}
\usepackage{cite}
\usepackage[toc,page]{appendix}
\usepackage{graphicx}
\usepackage[usenames,dvipsnames]{color}
\usepackage{authblk}

\usepackage[polutonikogreek,english]{babel}
\newtheorem{theorem}{Theorem}
\newtheorem{corollary}{Corollary}[theorem]

\begin{document}
\title{On the degree growth  of iterated birational maps}

\author{Claude M. Viallet\footnote{claude.viallet@upmc.fr} }
\affil{  Sorbonne Universit\'e,  Centre National de   la Recherche Scientifique}

\maketitle

\begin{abstract}
  We construct a family of birational maps acting on two dimensional
  projective varieties, for which the growth of the degrees of the
  iterates is cubic. It is known that this growth can be bounded,
  linear, quadratic or exponential for such maps acting on two
  dimensional compact K\"ahler varieties. The example we construct
  goes beyond this limitation, thanks to the presence of a singularity
  on the variety where the maps act. We provide all details of the
  calculations.
\end{abstract}

\section{Introduction}

A simple characterisation of the complexity of birational self-maps of
projective spaces is given by the nature of the growth of the sequence
of degrees $\{ d_n\}$ of their iterates: generic maps have exponential
growth, and some exceptional maps have polynomial
growth~\cite{Ve92,FaVi93,DiFa01}, \cite{Yo87,Ar90,RuSh97},  the latter
behaviour being linked with the notion of discrete
integrability~\cite{HiJoNi16}.

A straightforward test of integrability is then the vanishing of the
algebraic entropy~\cite{BeVi99} $$\epsilon = \lim_{n\rightarrow
  \infty}\frac{1}{n} \log( d_n).$$ This global index of complexity (and
equivalently its exponential called the dynamical degree) is always defined,
and is conjectured to be the logarithm of an algebraic integer for any
birational map~\cite{BeVi99} (see also~\cite{BeInJoMaSiTu19}). It is
clearly invariant by birational changes of coordinates. Moreover, if
the degree growth happens to be polynomial (vanishing entropy),
i.e. $d_n \simeq \alpha n^\nu$ when $n\rightarrow \infty$, then the
power $\nu$ is itself invariant by birational changes of coordinates.

The singularity structure is known to play an essential role in the
behaviour of the iterates.  Singularities appear in two instances:

-- birational self-maps of projective spaces of degree higher than $2$
have a non empty singular locus. A large number of results, most
notably the Lie algebraic characterisation of discrete Painlev\'e
equations~\cite{Sa01} inspired by \cite{Ok79}, have been obtained from
the description of this singular locus and the construction by blow
ups of rational varieties where the maps are
diffeomorphisms~\cite{Du10}.

-- one may also consider the action of maps on algebraic subvarieties
of projective spaces, and these varieties may or may not be smooth.

Both these types of singularities (of maps and of varieties) are
present in our analysis.

We construct a family of birational maps acting on two dimensional
projective varieties, having vanishing entropy, for which the growth
of the degrees of the iterates is {\em cubic}.  This goes beyond the
result of~\cite{DiFa01}, which would prevent anything between
quadratic and exponential growth. The existence of a singular point on
the varieties where the maps act is playing an important role in the
result. Such a behaviour has been observed
previously~\cite{AnMaVi02b,JoVi17,GuJoTrVi18}. Our point here is to
provide a complete proof that the growth is cubic, using a method
which can be applied in all dimensions and for non-vanishing entropy.

Our construction proceeds as follows: we start from a pencil of
elliptic curves in $\mathbb C \mathbb P_2$. All curves of this pencil
are invariant by a birational map $\varphi$, which has a specific
form, related to a second order recurrence relation. Such maps and
pencils are typical ingredients of the so called QRT
maps~\cite{QuRoTh88,QuRoTh89,Du10}. We then construct a third order
recurrence though the inflation process described
in~\cite{JoVi17}. This new recurrence defines birational maps in
$\mathbb C \mathbb P_3$, which automatically have a rational
invariant. The action on the two dimensional invariant projective
varieties obtained in this way is the one we analyse, proving that the
growth of the degree of the iterates is cubic.

\section{A pencil of elliptic curves}
Consider the recurrence relation
\begin{eqnarray}
  \label{order2}
  x_{n+1} + x_{n-1} = {\frac{a+b\,x_{n}+d\,{x_{n}}^{2}}{c+d\,x_{n}}}.
\end{eqnarray}
where $a,b,c,d$ are independent parameters.  It defines a birational
map on the projectivised space of initial conditions $\mathbb C
\mathbb P_2$:
\begin{eqnarray}
\label{phi}
  \varphi: [\xi, \zeta, \tau] \mapsto [-a\,{\tau}^{2}- \left(
    b\,\xi+c\,\zeta \right) \tau -d\; \xi\, \left( \xi+\zeta \right)
    ,\xi\,(c\,\tau+d\,{\xi}),\tau\,(c\,{\tau}+d\,\xi)]
\end{eqnarray}
written in homogeneous coordinates $[\xi, \zeta, \tau]$ of $\mathbb C
\mathbb P_2$. The map $\varphi$ has a rational invariant $\Pi$
\begin{eqnarray}
  \label{pi}
 {\tau}^{3}\; \Pi =   a \left( \xi+\zeta \right) {\tau}^{2}
  +b\,\tau\,\xi\,\zeta
  +c \left( {\xi}^{2} +{\zeta}^{2} \right) \tau
  +d\; \xi\,\zeta\, \left( \xi+\zeta \right) 
\end{eqnarray}
defining a pencil of invariant  elliptic curves in $\mathbb C \mathbb P_2$.

The map $\varphi$ is integrable, but we will not go into this
aspect. We will concentrate on the nature of the growth of the degrees
of the iterates of the maps we consider. For (\ref{phi}), this
sequence is
\begin{eqnarray}
  1,2,4,8,13,20,28,38,...
\end{eqnarray}
which has quadratic growth. The proof results from the construction of
an elliptic surface $R$ obtained by a finite number of blow ups of
$\mathbb C \mathbb P_2$, and the properties of the linear map induced
by $\varphi$ on the Picard group of $R$~\cite{Du10,Ta01,Mc07}.

\section{Inflation}
We apply to (\ref{order2}) the inflation process described in~\cite{JoVi17}
with the  simple transformation
\begin{eqnarray}
x_n = y_n + y_{n-1}
\end{eqnarray}
We get a recurrence relation of order three on $y_n$, and the
corresponding map in $\mathbb C \mathbb P_3$ is \\ $ \tilde{\varphi}:
[ x,y,z,t] \mapsto [x',y',z',t']$ with:
\begin{eqnarray} \label{phitilde}
  \begin{cases}
    & x'= -a\,{t}^{2}-b\, \left( x+y \right) t- c\, \left( x+y+z \right) t-d\, \left( x+    y \right)  \left( 2\,x+2\,y+z \right), \\
&     y'=x \,( c\, t+ d \left( x+y \right)), 
     z'= y\, ( c\, t+ d \left( x+y \right)), 
      t'= t\, ( c\, t+ d \left( x+y \right)).
  \end{cases}
\end{eqnarray}
which automatically has an invariant, constructed from the one
of $\varphi$ by  setting $ \xi=x+y$, $\zeta=y+z$, and  $\tau=t$.

In other words the new map $\tilde{\varphi}$ leaves the
pencil of two dimensional  surfaces $ \tilde{\Pi}(k) =0$  invariant
\begin{eqnarray}
  \label{pitilde}
 \tilde{\Pi}(k) = && k_0( a\,{t}^{2} \left( 2\,y+x+z \right)
  +  b\,  t \left( z+y \right)  \left( x+y \right) 
    +c \, t \left( (x+y)^{2} + ( y+z)^2\right)  \\ \nonumber
&&  +\, d \, \left( z+y \right)  \left( x+y \right)  \left( 2\,y+x+z \right))
  + { k_1} \, t^3.
\end{eqnarray}
with $k=(k_0,k_1)$

{\bf Lemma 1}
The varieties  $Q(k)$ of equation $ \tilde{\Pi}(k) =0$  are ruled surfaces. They are
cones with one singular point, the tip $\Sigma=[1,-1,1,0]$.

{\bf Proof:} Choosing a point $M \in \mathbb C \mathbb P_3$ (different
from $\Sigma$) determines a value $k(M)$ such that $M \in Q(k(M))$. It
is immediate to see that the line $\Lambda(M) = \tau \Sigma + (1-\tau)
M$ is entirely contained in $Q(k(M)))$. It is also straightforward to
check that the surfaces $Q(k)=0$ have $\Sigma$ as a unique singular
point.

{\bf Lemma 2}
The map $ \tilde{\varphi}$ sends  the line $  \Lambda(M)$ onto the line  $  \Lambda( \tilde{\varphi}(M))$

{\bf Proof:}
Straightforward.

Lemma 1 and 2 give the geometrical explanation of the fact that,
loosely speaking, $\Sigma$ is a fixed singular point of
$\tilde{\varphi}$.  We will give a precise meaning to this statement
in section \ref{fate}.

{\bf Remark} The singular locus of $ \tilde{\varphi}$ is the line
$\Delta$, which is intersection of the two hyperplanes $\{x+y=0\}$ and
$\{t=0\}$.  Notice that $\Delta \subset Q(k)$ for all $k$.

\section{The degree sequence, a fast analysis}
\label{fast}
The simplest way to get the exact value of the degree of the iterates
of $\tilde{\varphi}$ is to compute explicitly a sufficient number of
terms of the sequence of degrees, and try to fit the generating
function of their sequence with a rational function. We get
\begin{eqnarray}
  \label{deg_seq}
\{d_n\} =   1,2,4,8,14,23,35,51,71,96,126,162,204, \dots
\end{eqnarray}
which is fitted by the rational generating  function
\begin{eqnarray}
  \label{generating}
g(s) = \sum_{k=0}^{\infty} d_k s^k = {\frac {1-s +2\,s^3-{s}^{4}}{ \left( 1+s \right)  \left( 1-s \right) ^{4}}}
\end{eqnarray}
showing vanishing entropy~\cite{BeVi99} with {\em cubic } growth,
since all poles of $g(s)$ are on the unit circle and the
order of the pole at $s=1$ is 4.

It is straightforward to conjecture that the sequence of degrees is given by
\begin{eqnarray}
  d_n = {{17}/{16}}+{{5}/{12}} \,n +3/8\,{n}^{2}+1/12\,{n}^{3}-1/16\,
 \left( -1 \right) ^{n}
\end{eqnarray}
We will now prove this result.

\section{Preliminaries}
\label{prelim}

\subsection{The ${\bf K}$ varieties}

The rational map $\tilde{\varphi}$ has a rational inverse
$\tilde{\psi}$, which can be written $\tilde{\psi} = \lambda \cdot
\tilde{\varphi} \cdot \lambda$ with the linear involution $\lambda:
      [x,y,z,t] \mapsto [z,y,x,t]$.  The composition of
      $\tilde{\varphi} \cdot \tilde{\psi}$ is the identity, but when
      written with homogeneous coordinates, it appears as a
      multiplication of all the coordinates by a common factor:
\begin{eqnarray}
  \tilde{\psi} \cdot  \tilde{\varphi} ={\bf K}_{\tilde{\varphi}} \, id,
  \qquad
    \tilde{\varphi} \cdot  \tilde{\psi} = {\bf K}_{\tilde{\psi}} \, id.
\end{eqnarray}
The  factors  $ {\bf K}_{\tilde{\varphi}}$ and $  {\bf K}_{\tilde{\psi}}$ define varieties of codimension $1$  in  $\mathbb C \mathbb P_3$:
\begin{eqnarray}
  {\bf K}_{\tilde{\varphi}}= \left( c\,t+d \left( x+y \right)  \right) ^{3}  = B_1^3,
  \qquad
    {\bf K}_{\tilde{\psi}}=  \left( c\,t+d \left( z+y \right)  \right) ^{3} =C_1^3
\end{eqnarray}

We know that the map $ \tilde{\varphi}$ blows the hyperplane $\{B_1=0\}$ (resp. $
\tilde{\psi}$, $\{C_1=0\}$) down to something of smaller dimension and
the fate of the images under further iterates of $ \tilde{\varphi}$
(resp.  $ \tilde{\psi}$) eventually shapes the sequence of
degrees~\cite{BeVi99}.

The factors $ {\bf K}_{\tilde{\varphi}}$ also appear in the
calculations of the iterates (see section~\ref{two_ways}).

\subsection{Two ways of calculating iterates}

We will denote by $ \tilde{\varphi}^k $ the {\em regularised} $k$-th
power of $ \tilde{\varphi} $. By this we mean that when computing the
$k$th power with its homogeneous polynomial form, we remove all common
factors.

\label{two_ways}
We will make use of the seemingly trivial relation
  $ \tilde{\varphi}^{k+1}  =  \tilde{\varphi} \cdot \tilde{\varphi}^k  = \tilde{\varphi}^k \cdot \tilde{\varphi} $
to calculate iterates in two different ways:

1- Direct  image (i.e. $\tilde{\varphi}^{k+1} = \tilde{\varphi} \cdot \tilde{\varphi}^k$) \hfill\break
a) Apply the polynomial formulae giving the image by $\tilde{\varphi}$.  At
each step this will provide $N+1$ polynomials in the homogeneous
coordinates of the initial point. \hfill\break
b) Remove the common factors. \hfill\break
c) Proceed.

2 - Pull-back (i.e. $\tilde{\varphi}^{k+1} = \tilde{\varphi}^k \cdot \tilde{\varphi}$ ) \hfill\break
a) Get the  coordinates of $\tilde{\varphi}^{k+1} $  as the pull-backs by
$\tilde{\varphi}$ of the coordinates of $\tilde{\varphi}^{k} $. \hfill\break
b) Remove the common factors. \hfill\break
c) Proceed. 

\smallskip

Steps a) of the two methods give the same point in projective
space, but they possibly produce  different determinations, as  the
common factors differ in general.  Step b) is crucial. Once the common
factors are removed of course the determinations coincide.

The interest of the pull-back method  is the following:

Suppose $S$ is an indecomposable variety of codimension $1$ of
equation $E_S = 0$. The pullback by $\tilde{\varphi}$ of the equation
of $S$ gives the equation $E_{S'}$ of the image $S'$ of $S$ by
$\tilde{\psi}$.

\smallskip
{\bf Claim}: This pull-back may contain  additional
  factors, which are part of the { total} transform, but not part
of the { proper} transform, {\em  and the only possible
  factors are powers of $B_1$}.
\begin{eqnarray}
  \label{dual}
  \tilde{\varphi}^*( E_S )=  B_1^{\nu_S} \, E_{S'}
\end{eqnarray}
with $\nu_S$ some integer depending on $S$.

\section{ The fate of the ${\bf K}$  varieties under the action of  $ \tilde{\varphi}$}
\label{fate}

The map $ \tilde{\varphi}$ sends the whole hyperplane$\{ B_1 =0\}$ to
the point $ S_1= [1, 0, 0, 0]$. Further images are evaluated by first
a calculation of the image of a generic point and then the
specialisation of the coordinates to the ones of the running point
$[cx,cy,cz,- d \left( x+y \right)]$ of $B_1$. We get successively
$S_2= [2,- 1, 0, 0]$, $ S_3=[2, -2, 1, 0]$, $S_4= \Sigma = [1, -1, 1,
  0]$, and then $S_5=\Sigma$ again. The points $S_3$ and $\Sigma$ are
singular points of $ \tilde{\varphi}$. As mentioned above, $\Sigma$ is
a ``singular fixed point'' of $ \tilde{\varphi}$. By this we mean the
following:
\begin{theorem}
  \label{sigma}
  $\Sigma$ is a singular point of $ \tilde{\varphi}$  and $
  \tilde{\varphi}^{k+4} (B_1) = \Sigma$ for $k\geq 0$.
\end{theorem}
{\bf Proof of Theorem \ref{sigma}:} The proof uses blow-ups intended
to desingularise the map as in~\cite{CaTa18}. Since the hyperplane
$\{B_1=0\}$ is blown down to the point $S_1$ by $\tilde{\varphi}$, the
first step is to blow $S_1$ up and see what is the (new) image of $B_1$. To
be complete, the coordinatisation of the blow up should be done in
more than one coordinate patch but the main point is to ensure that
the induced map in the blown-up space does not send the two
dimensional hyperplane $\{B_1=0\}$ into a smaller dimensional space.

In the transformed space obtained by 
blowing up  $S_1$, the image of $\{B_1=0\}$ is a curve.  Therefore  we have to go on and blow-up again (along
that curve).  A possible choice of coordinates for the explicit
realisation of the blow ups is given in the
Appendix. The result of the calculation is that we need
three successive blow ups $\beta_{i,j}, j=1 \dots 3$ at each point
$S_i$. We denote by  $\delta_{i,j}$ the inverse of $\beta_{i,j}$ and by
$F_{i,j}$ the images of $\{B_1=0\}$ in the $j$-th blow up at $S_i$,
with evident notations (for example $F_{4,2}$ is the image of
$\{B_1=0\}$ by $ \tilde{\varphi}^4$ in $ \beta_{4,2} \cdot \beta_{4,1}
( \mathbb C \mathbb P_3)$).

We find that the $F_{i,1}$ and $F_{i,2}$ are curves, while the
$F_{i,3}$ are two dimensional.
   
To obtain the action of $\tilde{\varphi}^{k+4}$ on $\{B_1=0\}$
  for $k\geq 0$, we consider $\tilde{\varphi}\cdot
  \delta_{4,1}\cdot \delta_{4,2}\cdot \delta_{4,3} ( F_{4,3} )$. It is
  just $S_4=\Sigma$, which proves the theorem.

  {\bf Remark} The plane $\{w_{4,1}=0\}$ is sent to $\Sigma$ by
  $\tilde{\varphi}\cdot \delta_{4,1}$. The map on this plane induced
  by $\tilde{\varphi}$ is isomorphic to the one we started with
  (Eq. \ref{phi}).  One may indeed go further and blow up $S_5=S_4$
  with $\beta_{5,1} =\beta_{4,1}$. The image of $\{B_1=0\}$ by
  $\tilde{\varphi}^{5}$ is a curve $F_{5,1}$ (see Appendix) {\em This
    curve is different from $F_{4,1}$, of which it is the image by
    $\tilde{\varphi}$}.  As a consequence, if we wanted to
  desingularise the iterates of $\tilde{\varphi}$ by blow ups, one
  would need an infinite number of them, along the various curves
  $F_{n,1}$. Note that $F_{5,2}$ is a curve, and $F_{5,3}$ is two
  dimensional, as expected.

\begin{corollary} The action of the inverse $ \tilde{\psi}^k$ on the 
hyperplane $C_1$ being obtained via the linear symmetry $\lambda$,
$\Sigma$ is also a ``singular fixed point'' of $ \tilde{\psi}$ after
the $4$th iterate.
\end{corollary}


\section{Stabilisation of the iterates}
\label{stabil}
Start from a generic point $p_0=[x,y,z,t]$ and examine the
structure of the iterates of the map $\tilde{\varphi}$, calculated
explicitly, keeping track of the product structure of the components,
along the lines of~\cite{Vi15}, see also~\cite{KaMaTo18}:
\begin{eqnarray*}
 && p_1= [A_1, x\,B_1, y\,B_1, t\,B_1]\\
 && p_2= [A_2, A_1\,B_2, x\,B_1\,B_2, t\,B_1\,B_2]\\
 && p_3= [A_3, A_2\,B_3, A_1\,B_2\,B_3, t\,B_1\,B_2\,B_3]\\
 && p_4= [A_4, A_3\,B_4, A_2\,B_3\,B_4, t\,B_1\,B_2\,B_3\,B_4]\\
 && p_5= [A_5, A_4\,B_5, A_3\,B_4\,B_5, t\,B_1\,B_2\,B_3\,B_4\,B_5]
\end{eqnarray*}
    {\begin{theorem}
\label{theorem-stabil}
The general form of the iterates $p_k$ for $k\geq 3$ is
\begin{equation}  { p_k = [ A_k,\; A_{k-1}\,B_k,\;
      A_{k-2}\,B_{k-1}\,B_k,\; t\,B_1\,B_2\dots B_k]},
\end{equation}
or, defining $\quad \Gamma_k= t\,B_1\,B_2\dots B_k$, $\qquad p_k = [ A_k,\;
  A_{k-1}\,B_k,\; A_{k-2}\,B_{k-1}\,B_k,\; \Gamma_k]$.
\end{theorem}
    } We will prove Theorem~\ref{theorem-stabil} and give the
    recurrence relation between the successive $A_k$ and $B_k$ using
    both the direct way and the pull back way to express $p_{k+1}$ in
    terms of $p_k$. Denoting $[P^i_k], i=1..4$ the coordinates of
    $p_k$, we have
\begin{eqnarray}
 [P^i_{k+1}] \simeq \tilde{\varphi}([P^i_{k}]) \quad \mbox{(direct image)}\qquad and \qquad  [P^i_{k+1}] \simeq [\tilde{ \varphi}^*(P^i_k)])\quad \mbox{(pull-back)},
\end{eqnarray}
where the equivalence $ \simeq $ means equality up to  common factors.
\par\noindent
{\bf {Proof  of Theorem \ref{theorem-stabil}}:}
\par\noindent
Defining the sequences $\{\alpha\}$ and $\{\beta\}$ by
$ \tilde{\varphi}^*( A_k) = B_1^{\alpha(k)} \, A_{k+1} $ and 
$ \tilde{\varphi}^*( B_k) = B_1^{\beta(k)} \, B_{k+1} $,
we may evaluate  $\tilde{\varphi}^*(p_k)$ (up to common factors):
\begin{eqnarray}
  [{\bf B_1}^{\alpha(k)}\, A_{k+1},{\bf B_1}^{\alpha(k-1)+\beta(k)} \,
    A_{k}\,B_{k+1} ,\;{\bf B_1}^{\alpha(k-2)+\beta(k-1)+\beta(k)}\,
    A_{k-1}\,B_{k}\,B_{k+1},{\bf B_1}^{( \sum_{j=1}^k \beta(j))} \,
    \Gamma_{k+1} ]
\end{eqnarray}
We know from theorem \ref{sigma} that  no factor $B_1$ is left on
  the three first components for $k\geq 3$. This gives:
\begin{eqnarray}
  \alpha(k) = \alpha(k-1)+\beta(k)
  \end{eqnarray}
which in turn proves the stability of the form of $p_k$ since 
\begin{eqnarray}
 p_{k+1} = [ A_{k+1},\; A_{k}\,B_{k+1},\; A_{k-1}\,B_{k}\,B_{k+1},\;
   \Gamma_{k+1}].
\end{eqnarray}

\section{ Back to the sequence of degrees }

{\begin{theorem}
      \label{growth}
      The sequence of degrees $\{d_n\}$ is the one given in section
      \ref{fast}, i.e.
      \begin{eqnarray}
   d_n = {{17}/{16}}+{{5}/{12}} \,n +3/8\,{n}^{2}+1/12\,{n}^{3}-1/16\,
 \left( -1 \right) ^{n}
      \end{eqnarray}
\end{theorem}}

\noindent
{\bf Proof of Theorem \ref{growth}}:
We may  apply directly $\tilde{\varphi}$ to $p_k$, equate it
(up to factors) to $p_{k+1}$, and get the new recurrence relations
determining $A_{k+1}$ and $B_{k+1}$ in terms of the previous $A$'s and
$B$'s.
\begin{eqnarray}
  \begin{cases}
       & d \; D_k + c\;\Gamma_{k} - B_1 B_2 \dots B_{k-3} \, B_{k-2}^2
    \;{\bf B_{k+1}} =0 \\ & {\bf A_{k+1}} \Gamma_{k-2} B_{k-2} +
    a\,\Gamma_k^2 +B_kD_{k-1} ( d \,D_k + c\, \Gamma_k) + D_k \bigg(
    b\,\Gamma_k + d(D_k+A_k) \bigg) =0
  \end{cases}
\end{eqnarray}
with $D_k = A_k + A_{k-1}B_k$.

These relations give us the additional information on the degrees
$d_A$ and $d_B$ of the $A$'s and the $B$'s:
\begin{eqnarray*}
 d_B(n+1) - d_B(n) - d_B(n-1) + d_B(n-2) -1=0,\quad d_B(1) =1,\quad
 d_B(2)=2,\quad d_B(3)=4,
\end{eqnarray*}
giving
\begin{eqnarray*}
&d_B(n)& = 1/8+  n/2+1/4\,{n}^{2}-1/8\, \left( -1 \right) ^{n}, \\
&  d_A(n)& = 1 + \sum_{k=1}^n \, d_B(k),\\
&  d_A(n) & = {{17}/{16}}+{{5}/{12}} \,n +3/8\,{n}^{2}+1/12\,{n}^{3}-1/16\,
  \left( -1 \right) ^{n}.
\end{eqnarray*}

\textgreek{<'Oper >'edei de\^ixai}.

\section{Conclusion and perspectives}

\begin{itemize}
\item
  { The heuristic method used in section~\ref{fast} proves to be
    extremely efficient: one gets the results  in
  a fraction of a second with a small laptop.  The reason is the
  rational nature of the generating function of the sequence of
  degrees~(\ref{deg_seq}), a property which is not exceptional.
  Although it is not true for all birational maps (see~\cite{HaPr05}),
  it is verified for a large number of them. This moreover supports
  our conjecture that the entropy is the logarithm of an algebraic
  integer, {\em   keeping in mind that birationality is essential}, as is
  clear from sections \ref{prelim} and \ref{stabil}. We will address
  this question in a forthcoming publication.}
  \item
{ An important result of~\cite{DiFa01} is that the sequence of degrees
  of the iterates of birational maps on compact K\"ahler varieties of
  dimension 2 have bounded, linear, quadratic or exponential growth.
  What we have produced here is an example of birational maps acting on 
  2-dimensional algebraic hypersurfaces in $\mathbb CP_3$ and having cubic
  growth. The surfaces where the maps act have one singular point,
  and that point is a sink for the maps, allowing this different
  behaviour. One could also view this as an example of birational maps
  acting on a K\"ahler variety of dimension 3 with cubic growth.}
\item
{The fact that the growth is cubic implies in particular that the maps
  we constructed {\em cannot be obtained by a reduction of the known
    integrable quad equations\cite{AdBoSu03} since they all have at
    most quadratic growth~\cite{Vi08}.}}
\item
{Another interesting feature of the maps (\ref{phitilde}) is obtained
  by a simple graphical analysis of the orbits, showing the existence
  of an additional invariant, {\em which cannot be
    algebraic}~\cite{Gi80,Be99}.  It is worth noticing that the
  existence of a ``fixed singular point'' as shown in section
  \ref{fate} can be observed for the linearisable two dimensional
  example given in~\cite{AbHaHe00}. This similarity is not
  accidental. The integrability aspects will be treated elsewhere.}
\item
{The model can be made non-autonomous by setting $a=\alpha + \beta n +
  \gamma (-1)^n$. The sequence of degrees of the iterates {\em remains
    unchanged}, meaning that this non-autonomous generalisation is
  integrable as well.  The proof goes exactly along the same line as
  above.The fact that the sequence of degree is unchanged between the
  autonomous case and the non-autonomous case comes from the fact that
  the singularity structures are equivalent.}
\end{itemize}

{\bf{Acknowledgements}}

\bigskip
Parts of these results have been presented at the SIDE 13 meeting in
Fukuoka (Japan) November 2018 (http://www.side--conferences.net/).
The author would like to thank the Sydney Mathematical Research
Institute (University of Sydney), for hospitality and support during
the final elaboration of this work, and N. Joshi for a number of
challenging discussions and fruitful suggestions during the stay at
the SMRI.
\bigskip

{\bf{Appendix: Local coordinates for the blow-ups}}

\bigskip
 
  \label{appendA}
 We first have to deal with blow-ups at the points $S_i, i=1..4$.
 Denoting $[x,y,z,t]$ the original homogeneous components in $\mathbb
 C \mathbb P_3$, we see that they all have a non vanishing first
 coordinate $(x)$, and we will use local (affine) coordinates obtained
 by normalising this first coordinate to $1$.  There is a great
 arbitrariness in the choice of coordinates used to write the blow ups
 as explicit birational transformations.  A possible choice for the
 first blow-up is
 
    \begin{eqnarray}
      \beta_{1,1}:    [x,y,z,t] \mapsto [u_{1,1},v_{1,1},w_{1,1}]= [ z/y, t/y,t/x]
      \end{eqnarray}
    with  $[u_{1,1},v_{1,1},w_{1,1}]$ as  affine (local) coordinates  in the blown-up space.
    The image of $B_1$ is the curve  $F_{1,1} =\{ d\, ( u_{1,1}+1)+c \, v_{1,1}=0, \quad w_{1,1}=0\}$, and we need to perform a second  blow-up along this curve. This is done by
       \begin{eqnarray}
         \beta_{1,2}:   [u_{1,1},v_{1,1},w_{1,1}] \mapsto
         [u_{1,2},v_{1,2},w_{1,2}]=[{ u_{1,1}},{v_{1,1}},{\frac {{ w_{1,1}}}{d\, {u_{1,1}}+c\, v_{1,1}+ d }}]
       \end{eqnarray}
       The image of $B_1$ is still a curve $ F_{1,2} = \{ d\,( u_{1,2}+ 1) + c \,v_{1,2} =0,   ( ad -bc +c^2) \,  v_{1,2}  w_{1,2}+d=0  \}$,  and we need one more
         blow up to have $B_1$ sent to a two  dimensional image  by
         $\tilde{\varphi}$:
       \begin{eqnarray}
       \beta_{1,3}: [u_{1,2},v_{1,1},w_{1,2}] \mapsto [{ u_{1,2}},{
           v_{1,2}},\frac {( ad -bc +c^2) \,v_{1.2} \,  w_{1,2}+d}{ d\, u_{1,2}+c
             \,v_{1,2}+d  }]
       \end{eqnarray}
After this third  blow up,  $B_1$ is sent to   a two dimensional image $F_{1,3}$ by  $ \tilde{\varphi}$.
       \begin{eqnarray}
         F_{1, 3}  = \{  d\, ( u_{1,3} + 1) \, +c  \, v_{1,3}  =0 \}
     \end{eqnarray}   

       We may then proceed and blow up the point $S_2$. We again
       have to do three successive blow ups
       \begin{eqnarray}
         \begin{cases}
           \beta_{2,1}:[x,y,z,t] \mapsto [u_{2,1},v_{2,1},w_{2,1}]=[2\, z/(x+2\,y),2\, t/(x+2\,y),t/x]
         \\
     \beta_{2,2}:   [u_{2,1},v_{2,1},w_{2,1}]\mapsto  [u_{2,2},v_{2,2},w_{2,2}]= [ u_{2,1},{ w_{2,1}},{\frac { d( u_{2,1}+2 )+ ( b-2\,c )\, v_{2,1} }{{ w_{2,1}}}}]
         \\
         \beta_{2,3}: [u_{2,2},v_{2,2},w_{2,2}]\mapsto  [u_{2,3},v_{2,3},w_{2,3}]=
         [{ u_{2,2}},{ v_{2,2}},v_{2,2}/w_{2,2}]
         \end{cases}
       \end{eqnarray}

       Next are the blow ups over  $S_3$ and  $S_4$
  \begin{eqnarray}
         \begin{cases}
            \beta_{3,1}:[x,y,z,t] \mapsto [u_{3,1},v_{3,1},w_{3,1}]= [(2\,z-x)/2(x+y), t/(x+y), t/x]
         \\
     \beta_{3,2}:   [u_{3,1},v_{3,1},w_{3,1}]\mapsto  [u_{3,2},v_{3,2},w_{3,2}]=[u_{3,1}, v_{3,1}, w_{3,1} / ( v_{3,1}+d/c)]
         \\
         \beta_{3,3}:  [u_{3,2},v_{3,2},w_{3,2}]\mapsto  [u_{3,3},v_{3,3},w_{3,3}]= [u_{3,2},v_{3,2},\frac{   c\, v_{3,2}+d}{ 2 (ad-bc+c^2) w_{3,2} +c^2} ]
         \end{cases}
       \end{eqnarray}

  \begin{eqnarray}
         \begin{cases}
           \beta_{4,1}:[x,y,z,t] \mapsto [u_{4,1},v_{4,1},w_{4,1}]=[(z-x)/(x+y), t/(x+y),t/x]
         \\
     \beta_{4,2}:   [u_{4,1},v_{4,1},w_{4,1}]\mapsto  [u_{4,2},v_{4,2},w_{4,2}]=[{u_{4,1}},{v_{4,1}},{\frac { w_{4,1}}{d ( u_{4,1} +1)+c\, v_{4,1}}}]
         \\
         \beta_{4,3}:  [u_{4,2},v_{4,2},w_{4,2}]\mapsto  [u_{4,3},v_{4,3},w_{4,3}]= [ u_{4,2}, v_{4,2}, {\frac{  d( u_{4,2} +1)+ c v_{4,2} }{  2(ad-bc+c^2)  v_{4,2} w_{4,2}+d  }}]
         \end{cases}
 \end{eqnarray}
  \begin{eqnarray}
    \begin{cases}
   F_{4,1} =\{d\, u_{4,1} + c \, v_{4,1}+d =0 , \quad w_{4,1}   =0\}.
\\
  F_{4,2} =\{ d\, u_{4,2} + c \, v_{4,2}+d =0 ,  \quad
    v_{4,2}\, w_{4,2} +  2d/ (ad-bc+{c}^{2})   =0\}.
\\
  F_{4,3} =\{d\, u_{4,3} + c \, v_{4,3}+d =0 \}.
    \end{cases}
  \end{eqnarray}
  and
  \begin{eqnarray}
   F_{5,1} = \{ u_{5,1}^{2}d^{2}+d \left( b-c \right)  u_{5,1}\,v_{5,1}+3\,d
^{2} u_{5,1}+ \left( ad-{c}^{2} \right) v_{5,1}^{2}+bd \,v_{5,1}+ 2\,d^{2} , \quad w_{5,1}=0 \},
 \end{eqnarray}
 recalling that $[u_{4+k,1}, v_{4+k,1}, w_{4+k,1}] = [u_{4,1},
   v_{4,1}, w_{4,1}] $, but $[u_{4+k,j}, v_{4+k,j}, w_{4+k,j}] \neq
 [u_{4,j}, v_{4,j}, w_{4,j}]$ for $ j=2,3, k\geq0$.  The curve
 $F_{5,1}$ is rational.It is just the transform of $F_{4,1}$ by the
 map induced by $\tilde{\varphi}$.


\begin{thebibliography}{10}

\bibitem{Ve92}
A.P. Veselov, {\em Growth and Integrability in the Dynamics of Mappings}.
\newblock Comm.\ Math.\ Phys. {\bf 145} (1992), pp. 181--193.

\bibitem{FaVi93}
G.~Falqui and C.-M. Viallet, {\em Singularity, complexity, and
  quasi--integrability of rational mappings}.
\newblock Comm.\ Math.\ Phys. {\bf 154} (1993), pp. 111--125.
\newblock hep-th/9212105.

\bibitem{DiFa01}
J.~Diller and C.~Favre, {\em Dynamics of bimeromorphic maps of surfaces}.
\newblock Amer. J. Math. {\bf 123}(6) (2001), pp. 1135--1169.

\bibitem{Yo87}
Y.~Yomdin, {\em Volume growth and entropy}.
\newblock Israel J. Math. {\bf 57} (1987), pp. 285--299.

\bibitem{Ar90}
V.I. Arnold, {\em Dynamics of complexity of intersections}.
\newblock Bol. Soc. Bras. Mat. {\bf 21} (1990), pp. 1--10.

\bibitem{RuSh97}
A.~Russakovskii and B.~Shiffman, {\em Value distribution of sequences of
  rational mappings and complex dynamics}.
\newblock Indiana U. Math. J. {\bf 46} (1997), pp. 897--932.

\bibitem{HiJoNi16}
J.~Hietarinta, N.~Joshi, and F.W. Nijhoff.
\newblock {\em Discrete Systems and Integrability}.
\newblock Cambridge texts in applied mathematics. Cambridge University Press,
  (2016).

\bibitem{BeVi99}
M.P. Bellon and C-M. Viallet, {\em Algebraic Entropy}.
\newblock Comm. Math. Phys. {\bf 204} (1999), pp. 425--437.
\newblock chao-dyn/9805006.

\bibitem{BeInJoMaSiTu19}
R.~Benedetto, P.~Ingram, R.~Jones, M.~Manes, J.H.Silverman, and T.J. Tucker,
  {\em Current trends and open problem in arithmetic dynamics}.
\newblock Bull. Amer. Math. Soc.  (2019).
\newblock DOI: https://doi.org/10.1090/bull/1665.

\bibitem{Sa01}
H.~Sakai, {\em Rational Surfaces Associated with Affine Root Systems and
  Geometry of the {P}ainlev\'e Equations}.
\newblock Comm. Math. Phys. {\bf 220}(1) (2001), pp. 165--229.

\bibitem{Ok79}
K.~Okamoto, {\em Sur les feuilletages associ\'es aux \'equations du second
  ordre \`a points critiques fixes de P. {P}ainlev\'e}.
\newblock Jap. J. Math  (1979), pp. 1--79.

\bibitem{Du10}
J.J. Duistermaat.
\newblock {\em Discrete Integrable Systems: QRT Maps and Elliptic Surfaces}.
\newblock Springer Monographs in Mathematics. Springer New York,  (2010).

\bibitem{AnMaVi02b}
J-C.~Angl\`es d'Auriac, J-M. Maillard, and C-M. Viallet, {\em A classification
  of four-state spin edge {P}otts models}.
\newblock J. Phys. A: Math. Gen. {\bf 35} (2002), pp. 9251--9272.

\bibitem{JoVi17}
N.~Joshi and C-M. Viallet, {\em Rational maps with invariant surfaces}.
\newblock Journal of Integrable Systems {\bf 3}(1) (2018), p. xxy017.
\newblock arXiv:1706.00173.

\bibitem{GuJoTrVi18}
G.~Gubbiotti, N.~Joshi, D.T. Tran, and C-M. Viallet.
\newblock {\em Complexity and integrability in 4D bi-rational maps with two
  invariants}.
\newblock arXiv:1808.04942.

\bibitem{QuRoTh88}
G.R.W. Quispel, J.A.G. Roberts, and C.J. Thompson, {\em Integrable Mappings and
  Soliton Equations}.
\newblock Phys.\ Lett. {\bf A 126} (1988), p. 419.

\bibitem{QuRoTh89}
G.R.W. Quispel, J.A.G. Roberts, and C.J. Thompson, {\em Integrable Mappings and
  Soliton Equations II}.
\newblock Physica {\bf D34} (1989), pp. 183--192.

\bibitem{Ta01}
T.~Takenawa, {\em Discrete dynamical systems associated with root systems of
  indefinite type}.
\newblock Comm. Math. Phys. {\bf 224}(3) (2001), pp. 657--681.

\bibitem{Mc07}
C.T. McMullen, {\em Dynamics on blowups of the projective plane}.
\newblock Publ. Math. Inst. Hautes Etudes Sci. {\bf 105} (2007), pp. 49--89.

\bibitem{CaTa18}
A.S. Carstea and T.~Takenawa.
\newblock {\em Space of initial conditions and geometry of two 4-dimensional
  discrete Painlev\'e equations}.
\newblock arXiv:1810.01664.

\bibitem{Vi15}
C.~M. Viallet, {\em On the algebraic structure of rational discrete dynamical
  systems}.
\newblock J. Phys. A: Math. Theor. {\bf 48} (2015), p. 16FT01.

\bibitem{KaMaTo18}
M.~Kanki, T.~Mase, and T.~Tokihiro, {\em On the coprimeness property of
  discrete systems without the irreducibility condition}.
\newblock SIGMA {\bf 14} (2018), p. 065.

\bibitem{HaPr05}
B.~Hasselblatt and J.~Propp, {\em Degree-growth of monomial maps}.
\newblock Ergodic Theory and Dynamical Systems {\bf 27}(05) (2007), pp.
  1375--1397.
\newblock arXiv:math.DS/0604521.

\bibitem{AdBoSu03}
V.E. Adler, A.I. Bobenko, and Yu.B. Suris, {\em Classification of integrable
  equations on quad-graphs. The consistency approach}.
\newblock Comm. Math. Phys. {\bf 233}(3) (2003), pp. 513--543.
\newblock arXiv:nlin.SI/0202024.

\bibitem{Vi08}
C-M. Viallet, {\em Integrable lattice maps: $Q_V$, a rational version of
  $Q_4$}.
\newblock Glasgow Math. J. {\bf 51 A} (2009), pp. 157--163.
\newblock arXiv:0802.0294.

\bibitem{Gi80}
M.K. Gizatullin, {\em Rational ${G}$-surface}.
\newblock Izv. Akad. Nauk SSSR. Ser. Mat. {\bf 44}(1) (1980), pp. 110--144.
\newblock English translation is in "Math. USSR, Izvestiya", vol. 16 (1981),
  no.1, p.103-134: MR.

\bibitem{Be99}
M.P. Bellon, {\em Algebraic entropy of birational maps with invariant curves}.
\newblock Lett. Math. Phys. {\bf 50} (1999), pp. 79--90.

\bibitem{AbHaHe00}
M.J. Ablowitz, R.~Halburd, and B.~Herbst, {\em On the extension of the
  {P}ainlev\'e property to difference equations}.
\newblock Nonlinearity {\bf 13} (2000), pp. 889--905.

\end{thebibliography}
\end{document}